\documentclass[11pt,reqno]{article}
\usepackage[utf8]{inputenc}
\usepackage{booktabs} 
\usepackage{array} 
\usepackage{paralist} 
\usepackage{verbatim} 
\usepackage{amsthm}
\usepackage[table,xcdraw]{xcolor}
\usepackage[titletoc,toc,title]{appendix}
\usepackage{latexsym, amsmath, amssymb, a4, epsfig, color}
\usepackage{blindtext}
\usepackage{graphicx}
\usepackage{subcaption}
\usepackage[utf8]{inputenc}
\usepackage[export]{adjustbox}
\usepackage{wrapfig}
\usepackage{chngcntr}
\usepackage{apptools}
\usepackage{afterpage}

\usepackage{floatrow}

\usepackage{amssymb}
\usepackage{amsmath}
\usepackage[normalem]{ulem} 
\usepackage{algorithm}
\usepackage[noend]{algpseudocode}
\makeatletter
\def\BState{\State\hskip-\ALG@thistlm}
\makeatother
\AtAppendix{\counterwithin{definition}{subsection}}
\AtAppendix{\counterwithin{theorem}{subsection}}

\graphicspath{}

\newtheorem{theorem}{Theorem}[section]

\newtheorem{lemma}[theorem]{Lemma}

\newtheorem{definition}{Definition}

\newcommand{\Ga}{\alpha}
\newcommand{\Gb}{\beta}

\newcommand{\Gs}{\sigma}

\newcommand{\GO}{\Omega}

\newcommand{\RR}{\mathbb{R}}
\newcommand{\CC}{\mathbb{C}}
\newcommand{\NN}{\mathbb{N}}

\newcommand{\Om}{\Omega}
\newcommand{\ds}{\displaystyle}
\newcommand{\pf}{\noindent {\sl Proof}. \ }
\newcommand{\p}{\partial}
\newcommand{\pd}[2]{\frac {\p #1}{\p #2}}
\newcommand{\eqnref}[1]{(\ref {#1})}
\renewcommand{\qed}{\hfill $\Box$ \medskip}
\newcommand{\beq}{\begin{equation}}
\newcommand{\eeq}{\end{equation}}
\def\ep{\varepsilon}

\newcommand{\SingleOmega}{\mathcal{S}_{\partial\Omega}}

\newcommand{\Kcal}{\mathcal{K}}
\newcommand{\Scal}{\mathcal{S}}

\newcommand{\la}{\langle}
\newcommand{\ra}{\rangle}

\numberwithin{equation}{section}
\numberwithin{figure}{section}

\begin{document}

\newcommand{\TheTitle}{Shape reconstruction of a conductivity inclusion using the Faber polynomials}
\newcommand{\TheAuthors}{D. Choi, J. Kim and M. Lim}

\title{{\TheTitle}\thanks{{This work is supported by the Korean Ministry of Science, ICT and Future Planning through NRF grant No. 2016R1A2B4014530 (to D.C., J.K., and M.L.).}}}
\author{
Doosung Choi\thanks{\footnotesize Department of Mathematical Sciences, Korea Advanced Institute of Science and Technology, Daejeon 34141, Korea ({7john@kaist.ac.kr}, {kjb2474@kaist.ac.kr},  {mklim@kaist.ac.kr}).}\and Junbeom Kim \footnotemark[2] \and Mikyoung Lim\footnotemark[2]}

\date{\today}
\maketitle

\begin{abstract}

We consider the shape reconstruction of a conductivity inclusion in two dimensions. We use the concept of Faber polynomials Polarization Tensors (FPTs) introduced in \cite{choi:2018:GME} to derive an exact shape recovery formula for an inclusion with the extreme conductivity. This shape can be a good initial guess in the shape recovery optimization for an inclusion with either small or large conductivity values. We illustrate and validate our results with numerical examples.
\end{abstract}



\section{Introduction}

Let $\Omega$ in $\RR^2$ be simply connected, bounded, and Lipschitz. The background is homogeneous with conductivity $1$ and $\Omega$ is occupied with a material of dielectric constant $\sigma_0$ with $0 \le \sigma_0 \neq 1  \le \infty$. 
We consider
the following transmission problem of the quasi-static approximation of electromagnetic fields:
\beq\label{cond_eqn0}
\begin{cases}
\ds\nabla\cdot\sigma\nabla u=0\quad&\mbox{in }\RR^2, \\
\ds u(x) - H(x)  =O({|x|^{-1}})\quad&\mbox{as } |x| \to \infty
\end{cases}
\end{equation}
with $
\sigma=\sigma_0\chi(\Om)+\chi(\RR^2\setminus \overline{\Om})$,
where $H$ is an entire harmonic function and 
the symbol $\chi(D)$ indicates the characteristic function for a domain $D$. The solution $u$ should satisfy the transmission condition
$$u\big|^+=u\big|^-\quad\mbox{and }\quad \pd{u}{\nu}\Big|^+=\sigma_0\pd{u}{\nu}\Big|^-\qquad\mbox{on }\p\Om.$$
Here, $\nu$ is the outward unit normal vector on $\p\Om$ and the symbols $+$ and $-$ indicate the limit from the interior and exterior of $\Om$, respectively.

In this paper, we consider the problem of reconstructing the shape of $\Om$ from the measurements of $u$ away from $\Om$. Regarding the shape reconstruction, we refer readers to see previous papers in \cite{ammari:2012:GPT, ammari:2014:GPT}. We use the solution expression of $u$ associated with the exterior conformal mapping developed in \cite{jung:2018:NSS} and the geometric multipole expansion introduced in \cite{choi:2018:GME}. To explain the details, let's denote $\Psi$ the exterior conformal mapping associated with $\Om$. Indeed, according to the Riemann mapping theorem, there exist $\gamma>0$ uniquely and the conformal mapping $\Psi$ from $\{w\in\CC:|w|>\gamma\}$ onto $\CC\setminus\overline{\Om}$ such that
\beq\label{conformal:Psi}
\Psi(w)=w+a_0+\frac{a_1}{w}+\frac{a_2}{w^2}+\cdots.
\eeq
Here and after, we identify $z=x_1+ix_2$ in $\CC$ with $x=(x_1,x_2)$ in $\RR^2$. 
The conformal mapping coefficients $a_n$ can be solved numerically by the boundary integral equation involving the Neumann-Poincar\'{e} (NP) operator; see \cite{jung:2018:NSS}.

As a univalent function, the exterior conformal mapping $\Psi$ defines the so-called Faber polynomials, which is a basis for holomorphic functions in the domain associated with the exterior conformal mapping. We use the concept of Faber polynomials Polarization Tensors (FPTs) introduced in \cite{choi:2018:GME} to derive an exact shape recovery formula for an inclusion with the extreme conductivity. This shape can be a good initial guess in the shape recovery optimization for inclusions with either small or large conductivity values. We illustrate and validate our results with numerical examples.

The rest of the paper is organized as follows. In section 2, we formulate the transmission problem \eqnref{cond_eqn0} by boundary integrals and review the multipole expansions. Section 3 derives the exact inversion formula for the inclusion with the extreme conductivity. In section 4, we propose an optimization scheme with the initial guess that obtained from the exact shape recovery formula for the inclusion with the extreme conductivity. 

\section{Faber polynomial Polarization Tensors (FPTs)}

In this section we provide the boundary integral formulation for the transmission problem. We then explain the classical multipole expansion. Following \cite{choi:2018:GME}, we finally define the Faber polynomial Polarization Tensors (FPTs) as follows we define the Faber polynomial Polarization Tensors (FPTs) and the geometric multipole.


	Let $\Om$ be a simply connected and Lipschitz domain in $\RR^2$. 
	The single layer potential $\Scal_{\p \Om}$ and the Neumann-Poincar\'{e} (NP) operator $\Kcal_{\p \Om}^*$ associated with $\Om$ are defined as follows: for $\varphi\in L^2(\p\Om)$,
	\begin{align*}
	\ds&\Scal_{\p\Om}[\varphi](x)=\int_{\p D}\Gamma(x-\tilde{x})\varphi(\tilde{x})\, d\sigma(\tilde{x}),\quad x\in\RR^d,\\
	\ds &\Kcal_{\p\Om}^*[\varphi](x)=p.v.\,\frac{1}{2\pi}\int_{\partial\Omega}\frac{\left\la x-\tilde{x},\nu_x\right\ra}{|x-\tilde{x}|^2}\varphi(\tilde{x})\, d\sigma(\tilde{x}),\quad x\in\p \Om.
	\end{align*}
	Here, $\nu$ is the outward unit normal vector to $\p \Om$ and  $\Gamma(x)$ is the fundamental solution to the Laplacian, i.e., 
	$\Gamma(x)=(2\pi)^{-1}\ln|x|$.
	We also denote $\SingleOmega[\varphi](z):=\SingleOmega[\varphi](x)$ for $x=(x_1,x_2)$ and $z=x_1+ix_2$.
	
	One can express the solution to \eqnref{cond_eqn0} as
	\beq\label{umh1}
	u(x)=H(x)+\Scal_{\p\Om}[\varphi](x),\quad x\in\RR^2,
	\eeq
	where
	\beq\label{umh2}
	\varphi=(\lambda I-\Kcal_{\p\Om}^*)^{-1}\left[\nu\cdot \nabla H\right]\quad\mbox{with }\lambda = \frac{\sigma_0+1}{2(\sigma_0-1)}.
	\eeq
The invertibility of the operator $(\lambda I-\Kcal_{\p\Om}^*)$ is well-known for $|\lambda|\geq1/2$ as shown in \cite{escauriaza:1992:RTW,kellogg:2012:FPT,verchota:1984:LPR} (see also \cite{escauriaza:1992:RTW, escauriaza:1993:RPS}). We commend readers to see \cite{helsing:2013:SIE, helsing:2017:CSN} for the Numerical method to solve the integral equation and \cite{ammari:2004:RSI, ammari:2007:PMT} and references therein for more about the NP operator.

\subsection{Generalized Polarization Tensors (GPTs)}

	By applying the Taylor series expansion to \eqnref{umh1}, one can derive a multipole expansion for the transmission problem. 
	In terms of the conventional multi-index notation,
	$$x^{\alpha}=x_1^{\alpha_1}x_2^{\alpha_2},\quad
	|\alpha|=\alpha_1+\alpha_2,$$ the fundamental solution to the Laplacian and the background potential admit the Taylor series expansion
	\begin{align}
	\label{Taylor1}\Gamma(x-y)&=\sum_{|\alpha|=0}^\infty\frac{(-1)^{|\alpha|}}{\alpha!}\partial^\alpha \Gamma(x)y^\alpha,\\\label{Taylor2}
	H(y)&=\sum_{|\beta|=0}^\infty\frac{1}{\beta!}\partial^\beta H(0)y^\beta
	\end{align}
	for $y\in\p\Om$ and sufficiently large $x$. The Generalized Polarization Tensors (GPTs) associated with a domain $\Om$ and the conductivity $\sigma_0$ are defined as \beq
	\label{gpt}
	M_{\alpha\beta}(\Om,\lambda) = \int_{\p \GO} y^\alpha \left(\lambda
	I - \Kcal^*_{\p\GO}\right)^{-1}\left[\nu \cdot \nabla y^\beta\right](y) \, d\Gs (y)
\end{equation}
for the multi-indices $\alpha$, $\beta$.
By applying \eqnref{Taylor1} and \eqnref{Taylor2} into \eqnref{umh1} and \eqnref{umh2}, the multipole expansion for the solution to
\eqnref{cond_eqn0} is (see \cite{ammari:2007:PMT} for detail)
\beq\label{CP2}
u(x) = H(x)+\sum_{|\Ga|,|\Gb|=1}^\infty\frac{(-1)^{|\Ga|}}{\Ga!\Gb!}\p^\Ga\Gamma(x)M_{\Ga\Gb}(\Om,\sigma_0)\p^\Gb H(0), \quad |x|\gg1,
\eeq

\begin{definition}
Following \cite{ammari:2013:MSM}, we denote $\NN_{mk}^{(1)}$ and $\NN_{mk}^{(2)}$ as complex contracted GPTs:
\begin{align*}
\NN_{mk}^{(1)}(\Om, \lambda)&=\int_{\p\Om}z^k(\lambda I-\Kcal^*_{\p\Om})^{-1}\left[\pd{z^m}{\nu}\right]\,d\sigma(z),\\
\NN_{mk}^{(2)}(\Om, \lambda)&=\int_{\p\Om}z^k(\lambda I-\Kcal^*_{\p\Om})^{-1}\left[\pd{\overline{z^m}}{\nu}\right]\,d\sigma(z)\quad\mbox{for }m,k\in\NN.
\end{align*}
Here, $z=x_1+ix_2$ for $(x_1,x_2)\in \p \Om$. As the complex polynomials are linear combinations of real polynomials, the complex contracted GPTs are linear combinations of $M_{\alpha\beta}$'s.
\end{definition}

\subsection{Faber polynomials and FPTs}\label{sec:GMT}

As a univalent function, the exterior conformal mapping $\Psi$ defines the so-called Faber polynomials, $F_m(z)$'s, which is complex monomials and form a basis for complex analytic functions in $\Om$ (see \cite{duren:1983:UF}). 
The Faber polynomials are first introduced by G. Faber in \cite{faber:1903:UPE} and have been extensively studied in various areas.

The Faber polynomials $\{F_m(z)\}$ associated with $\Psi$ are defined by the relation
	\beq\label{eqn:Fabergenerating}
	\frac{w\Psi'(\zeta)}{\Psi(w)-z}=\sum_{m=0}^\infty \frac{F_m(z)}{w^{m}},\quad z\in{\overline{\Om}},\ |w|>\gamma.
	\eeq
The complex logarithm admits the series expansion (see \cite{duren:1983:UF, faber:1903:UPE, jung:2018:NSS}): for ${z}=\Psi(w)\in\CC\setminus\overline{\Om}$ and $\tilde{z}\in\Om$,
\beq\label{log:Faber}
\log({z}-\tilde{z})=\log w-\sum_{m=1}^\infty \frac{1}{m}F_m(\tilde{z})w^{-m}\eeq
with a proper branch cut. 
Each $F_m$ is an $m$-th order monic polynomial. 
	For example, the first three polynomials are
	$$F_0(z)=1,\quad F_1(z)=z-a_0,\quad F_2(z)=z^2-2a_0 z+a_0^2-2a_1.$$
	In general, if denote the Faber polynomial as $F_m(z) = \sum_{n=0}^{m} a_{mn} z^n$, we can find $\{a_{mn}\}$ recursively by the relation
\beq\label{Faberrecursion}
F_{n+1} (z) = z F_n (z) - n a_n - \sum_{s=0} ^{n} a_s F_{n-s} (z), \quad n\ge 0.
\eeq
By substituting $z=\Psi(w)$, we obtain
	\begin{equation} \label{eqn:Faberdefinition}
	F_m(\Psi(w))
	=w^m+\sum_{k=1}^{\infty}c_{mk}{w^{-k}},
	\end{equation}
	where $c_{mk}$'s are so-called Grunsky coefficients. 
In view of the Laurent series expansion \eqnref{conformal:Psi}, one can easily see by setting $m=1$ that
$$c_{1k}=a_k\quad\mbox{for all }k\in\NN.$$ 
	Recursive relation for general $m$ is also well-known (see \cite{duren:1983:UF}).

\begin{definition}
Following \cite{choi:2018:GME}, we define the Faber polynomial Polarization Tensors  (FPTs) as follows:
	For $m,k\in \NN$, we define
	\begin{align}\label{eqn:FPT1}
	F_{mk}^{(1)}(\Om, \lambda)=\int_{\p\Om}F_k(z)\left(\lambda I-\Kcal^*_{\p\Om}\right)^{-1}\left[\pd{F_m}{\nu}\right](z)\,d\sigma(z),\\\label{eqn:FPT2}
	F_{mk}^{(2)}(\Om, \lambda)=\int_{\p\Om}F_k(z)\left(\lambda I-\Kcal^*_{\p\Om}\right)^{-1}\left[\pd{\overline{F_m}}{\nu}\right](z)\,d\sigma(z).
	\end{align}
	We call $F_{mk}^{(1)}(\Om,\lambda)$ and $F_{mk}^{(2)}(\Om,\lambda)$ the Faber polynomial Polarization Tensors (FPTs) associated with $\Om$.
\end{definition}

By using the expansion of the complex logarithm \eqnref{log:Faber} and the integral formulation \eqnref{umh1}, one can express the solution to the transmission problem \eqnref{cond_eqn0} as follows  \cite{choi:2018:GME}:
	for a given harmonic function $H(z) = \sum_{m=1}^\infty \left(\alpha_m F_m(z)+\beta_m\overline{F_m(z)}\right)$ with complex coefficients $\alpha_m$ and $\beta_m$, the solution $u$ to \eqnref{cond_eqn0} satisfies that for $z=\Psi(w)\in\CC\setminus\overline{\Om}$,
	\begin{align}\label{eqn:FPT_expan}
	u(z)=H(z)-\sum_{k=1}^\infty \sum_{m=1}^\infty\frac{1}{4\pi k}\left[ \left( \alpha_m F_{mk}^{(1)}+\beta_m  \overline{F_{mk}^{(2)}}
	\right){w^{-k}}
	+ \left(\alpha_m  F_{mk}^{(2)}+\beta_m \overline{F_{mk}^{(1)}}\right)\overline{w^{-k}}\right].
	\end{align}
It is worth highlighting that the geometric expansion holds for all $z$ in the exterior of $\Om$ while the classical multipole expansion holds for $|z|\gg 1$.

The FPTs contains the information on the material parameter and the shape of the inclusion $\Om$. Actually, one can express the FPTs in terms of the Grunsky coefficients of $\Om$ and $\lambda$ (or, in terms of $\sigma_0$).
\begin{lemma}[\cite{choi:2018:GME}]\label{lemma:FPT}
Let $C$ be the Grunsky matrix with elements $\{c_{mk}\}$, and $\gamma^{\pm2\mathbb{N}}$ be diagonal matrices whose $(k,k)$-entry is $\gamma^{\pm2k}$.
For each $m,k\in\NN$, 
the FPTs satisfy
\begin{align*}
F_{mk}^{(1)} (\Om, \lambda)&=  4\pi k c_{mk} + 4\pi k \left( \frac{1}{4} -\lambda^2 \right) \left[ \gamma^{2\mathbb{N}} \left( \lambda^2 I - \Bigl(\frac{\gamma^{-2\mathbb{N}}C}{2}\Bigr)^2 \right)^{-1}\gamma^{-2\mathbb{N}}C\right]_{mk},\\
F_{mk}^{(2)}(\Om, \lambda) &= 8\pi k\lambda \gamma^{2k}\delta_{mk} + 4\pi k \left( \frac{1}{4} -\lambda^2 \right) \left[2\lambda\gamma^{2\mathbb{N}} \left( \lambda^2 I - \Bigl(\frac{\gamma^{-2\mathbb{N}}C}{2}\Bigr)^2 \right)^{-1} \right]_{mk}.
\end{align*}
Here, $\delta_{mk}$ is the Kronecker delta function.
\end{lemma}


\section{Exact inversion formula from the GPTs of extreme conductivities}
In this section, we provide the exact inversion formula for the extreme conductivities $\sigma_0=\infty \mbox{ or } 0$ (in other words, $\lambda = \pm 1/2$) from the GPTs. We derive the formula by using Lemma \ref{lemma:FPT}.

\begin{theorem}\label{thm:exact}
We have the following formulas.
\begin{itemize}
\item[\rm(a)]
For $m \ge 1$,
\begin{align*}
\gamma = \sqrt{\pm \frac{F_{11}^{(2)}(D,\pm \tfrac{1}{2})}{4\pi}}, \quad a_m = \frac{F_{m1}^{(1)}(D,\pm \tfrac{1}{2})}{4\pi m}.
\end{align*}

\item[\rm(b)]
For the extreme conductivity cases, i.e., when $\lambda = \pm\frac{1}{2}$, the coefficients of the exterior conformal mapping can be derived as follows.
\beq\label{recur1}
\gamma = \sqrt{\pm \frac{\mathbb{N}_{11}^{(2)}(D,\pm \frac{1}{2})}{4\pi}},  \quad a_0 = \frac{\mathbb{N}_{21} ^{(2)} (D,\pm \frac{1}{2}) }{2 \mathbb{N}_{11} ^{(2)} (D,\pm \frac{1}{2}) },  \quad a_1 = \frac{\mathbb{N}_{11}^{(1)} (D,\pm \frac{1}{2}) }{4\pi}.
\eeq
For $m \ge 2$, 
\beq a_m = \frac{\sum_{n=1} ^{m} a_{mn} \mathbb{N}_{n1} ^{(1)} (D,\pm \frac{1}{2})}{4\pi m},\label{recur2}
\eeq
where we get $\{a_{mn}\}_{n=1}^m$ from \eqnref{Faberrecursion}.
\end{itemize}
\end{theorem}

\pf
From Lemma \ref{lemma:FPT}, when $\lambda = \pm \frac{1}{2}$, the FPTs satisfy
	\begin{align*}
	F_{m1}^{(1)}(D,\pm \tfrac{1}{2})&
	=4\pi c_{m1} = 4\pi m a_m \qquad\textrm{ for } m \ge 1,\\
	F_{11}^{(2)}(D,\pm \tfrac{1}{2})& =\pm 4\pi \gamma^{2},\quad F_{m1}^{(2)}(D,\pm \tfrac{1}{2}) = 0 \qquad\textrm{ for } m \ge 2.
	\end{align*}
Recall that $F_m(z) = \sum_{n=0}^{m} a_{mn} z^n$. Note that $a_{mm}=1$ and the coefficients $\{a_{mn}\}_{n=1}^m$ can be induced by $\{ a_0, a_1, \cdots , a_{m-1}\}$ using \eqnref{Faberrecursion}. From the definition, the FPTs can be expressed as
	\begin{align*}
	F_{m1}^{(1)}(D,\lambda) = \sum_{n=1} ^{m} a_{mn} \mathbb{N}_{n1} ^{(1)} (D,\lambda), \quad F_{m1}^{(2)}(D,\lambda) = \sum_{n=1} ^{m} a_{mn} \mathbb{N}_{n1} ^{(2)} (D,\lambda).
	\end{align*}
	Hence, we obtain $\gamma$, $a_0$, and $a_1$ as follows.
	\begin{align}
	\left. \begin{array}{l}
	\ds \gamma =\sqrt{\pm \frac{F_{11}^{(2)}(D,\pm \frac{1}{2})}{4\pi}} =\sqrt{\pm \frac{\mathbb{N}_{11}^{(2)}(D,\pm \frac{1}{2})}{4\pi}}, \\[2mm]
	\ds a_0 = \frac{\mathbb{N}_{21} ^{(2)} (D,\pm \frac{1}{2}) }{2 \mathbb{N}_{11} ^{(2)} (D,\pm \frac{1}{2}) },\\[4mm]
	\ds a_1 = \frac{ F_{11}^{(1)}(D,\pm \frac{1}{2})}{4\pi} = \frac{\mathbb{N}_{11}^{(1)} (D,\pm \frac{1}{2}) }{4\pi}. \end{array} \right.
	\end{align}
For $m \ge 2$, each $a_m$ is derived from the relation $F_{m1} ^{(1)} (D,\pm \frac{1}{2}) = 4\pi m a_m$.
\qed

For the insulating inclusion, a recursion formula similiar to \eqnref{recur1} and \eqnref{recur2} was derived in \cite{choi:2018:CEP}.

\begin{theorem}
	Suppose that the inclusion has the extreme conductivity, i.e., $\lambda = \frac{1}{2}$ or $-\frac{1}{2}$. For given GPTs, $\widetilde{M}_{\alpha\beta}$ with $|\alpha|, |\beta| \le N$, there exist a unique simply connected domain $D$ that satisfies
	\begin{itemize}
		\item[(i)] $M_{\alpha \beta} (D, \lambda) = \widetilde{M}_{\alpha\beta}$ for $|\alpha|, |\beta| \le N$.
		\item[(ii)] The exterior conformal mapping $\Psi$ of $D$ has finite terms of order at most $N$, {\it i.e.}
		$$\Psi(w) = w+ a_0 + \frac{a_1}{w} + \cdots + \frac{a_N}{w^N}.$$
	\end{itemize}
\end{theorem}

\begin{figure}[H]	
	\begin{subfigure}[b]{0.3\textwidth}
		\includegraphics[width=\textwidth,trim=40 20 40 0, clip]{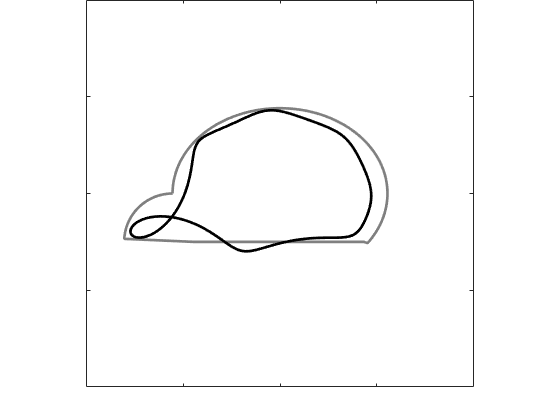}
	\end{subfigure}
	\begin{subfigure}[b]{0.3\textwidth}
		\includegraphics[width=\textwidth,trim=40 20 40 0, clip]{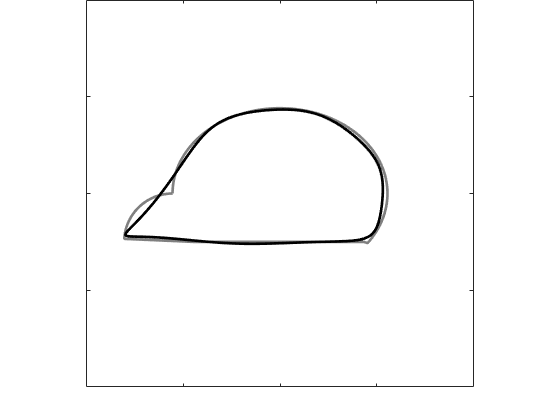}
	\end{subfigure}
	\begin{subfigure}[b]{0.3\textwidth}
		\includegraphics[width=\textwidth,trim=40 20 40 0, clip]{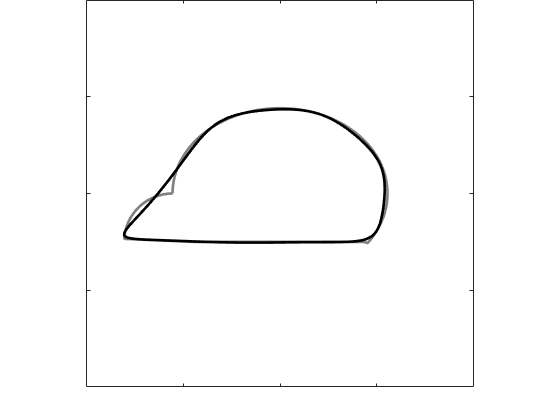}
	\end{subfigure}
	\caption{The initial  guess for the cap-shaped domain. The black curve is the reconstructed curve and the gray curve is the actual shape. Starting from the left side to the right side, the conductivity is 10, 50, 100.}
	\label{fig:cap_initial}
\end{figure}

\begin{figure}[H]	
	\begin{subfigure}[b]{0.3\textwidth}
		\includegraphics[width=\textwidth,trim=40 20 40 0, clip]{cap_confmap_condD50}
	\end{subfigure}
	\begin{subfigure}[b]{0.3\textwidth}
		\includegraphics[width=\textwidth,trim=40 20 40 0, clip]{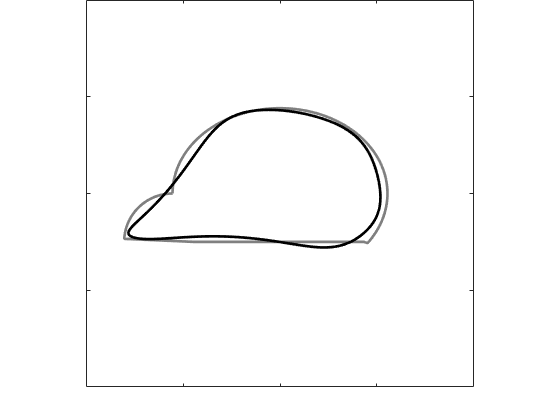}
	\end{subfigure}
	\begin{subfigure}[b]{0.3\textwidth}
		\includegraphics[width=\textwidth,trim=40 20 40 0, clip]{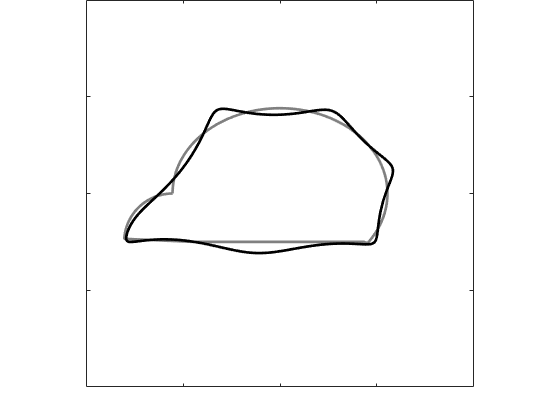}
	\end{subfigure}
	\caption{The initial  guess with errors. The black curve is the reconstructed curve and the gray curve is the actual shape. Starting from the left side to the right side, the results are reconstructed shapes from the data with 0\%, 10\%, 20\% error. Here, the conductivity $k= 50$.}
	\label{fig:initial_error}
\end{figure}

Below Figure \ref{fig:cap_initial} shows the shape recovered from the exact formula in Theorem \ref{thm:exact} for various $\sigma_0$. When $\sigma_0$ is large, the shape recovery is very close to the true shape. 
The example in Figure \ref{fig:cap_initial} is the initial guess with extreme conductivity for cap-shaped domain. Even the target domain has corners, the initial guess with extreme conductivity is quite accurate.  Figure \ref{fig:initial_error} reveals the initial guess with errors. For each GPT, we put the random error at most 0\%, 10\%, 20\%. The results in Figure \ref{fig:initial_error} shows the stability of the method.



\section{Optimization scheme to reconstruct shape from the GPTs}\label{sec:shape}

The exact recovery formula obtined in the previous section holds only for the extreme case $\lambda=\pm 1/2$. 
However, we may use the shape from the formula, assuming $\lambda=\pm 1/2$, as an initial guess when $\lambda$ is close to $\pm 1/2$. We call such a shape the reference shape for the inclusion. 
In the following we compare the reference shape with the equivalent ellipse, which is commonly used as the initial guess for the shape recovery problem, and then provide the optimization method to recover the shape for general conductivity value.

\subsection{Initial guess}\label{sec:initial}
\noindent\textbf{Equivalent ellipse}
One way to make initial guess is to find the equivalent ellipse, which is explained in \cite{ammari:2004:RSI}. When the first order GPTs $$M = \left[ \begin{array}{cc} m_{11} & m_{12} \\ m_{21} & m_{22} \end{array} \right]$$ are given, we can compute the equivalent ellipse. Suppose that the eigenvalues of $M$ are $\lambda_1, \lambda_2$ with $\lambda_1 > \lambda_2$ and corresponding eigenvectors are $(e_{11}, e_{12})^{T}$ and $(e_{21}, e_{22})^{T}$. Then for $\frac{1}{p} = \frac{\sigma_0-1}{\sigma_0+1} \left( \frac{1}{\lambda_1} + \frac{1}{\lambda_2} \right)$ and $q=\frac{\lambda_2 - \sigma_0 \lambda_1}{\lambda_1 - \sigma_0 \lambda_2}$,
$$a= \sqrt{\frac{p}{\pi q}},\quad b= \sqrt{\frac{pq}{\pi}},\quad \theta = \arctan \left( \frac{e_{21}}{e_{11}} \right).$$

The location of inclusion can be approximated from the second order GPTs since the first order GPTs are invariant under translation.
\smallskip

\noindent\textbf{Reference shape}
In this paper, we suggest new initial guess based on Theorem \ref{thm:exact}. From Lemma \ref{lemma:FPT}, we can express the conformal mapping coefficient as follows,
$$a_m = \frac{ F_{m1}^{(1)}(B,\pm \frac{1}{2})}{4\pi} + O\left(|\lambda|-\frac{1}{2}\right).$$
Hence, for nearly extreme conductivity case, let $\lambda = \frac{1}{2}$ when $\sigma_0 \gg 1$ and $\lambda = -\frac{1}{2}$ when $\sigma_0 \ll 1$. Then calculating $\gamma$ and conformal mapping coefficients with the method on the proof of Theorem \ref{thm:exact} become a good initial guess.

\begin{figure}[H]
\begin{subfigure}[b]{0.3\textwidth}
\includegraphics[width=\textwidth,trim=40 20 40 0, clip]{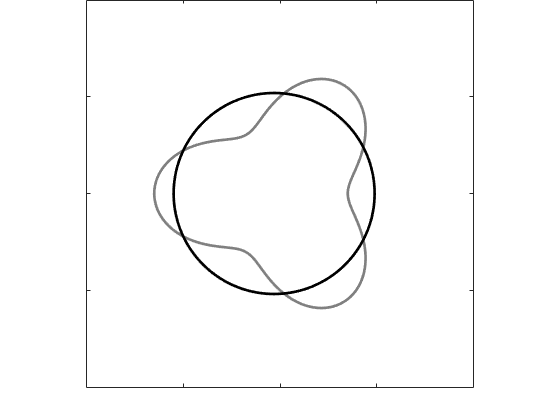}
\end{subfigure}
\begin{subfigure}[b]{0.3\textwidth}
\includegraphics[width=\textwidth,trim=40 20 40 0, clip]{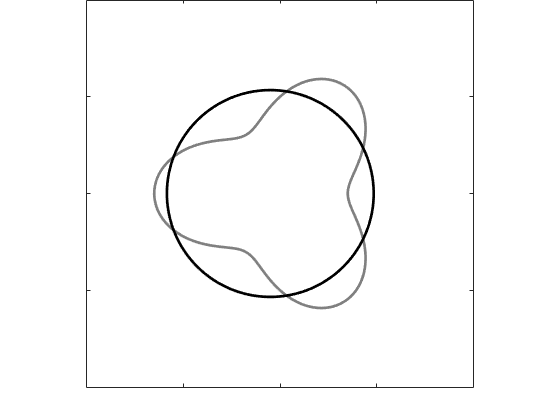}
\end{subfigure}
\begin{subfigure}[b]{0.3\textwidth}
\includegraphics[width=\textwidth,trim=40 20 40 0, clip]{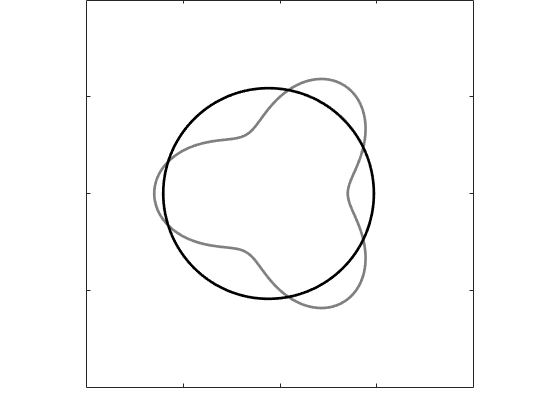}
\end{subfigure}
\vskip 0.2cm
\begin{subfigure}[b]{0.3\textwidth}
\includegraphics[width=\textwidth,trim=40 20 40 0, clip]{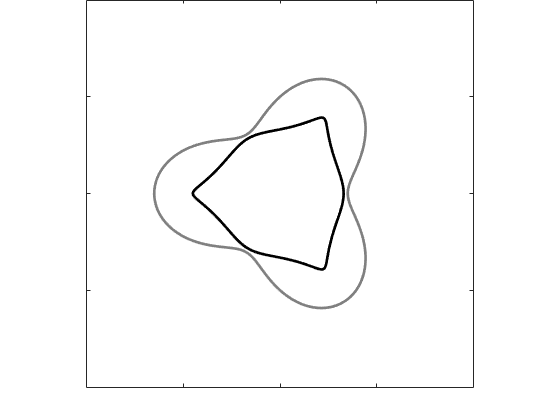}
\end{subfigure}
\begin{subfigure}[b]{0.3\textwidth}
\includegraphics[width=\textwidth,trim=40 20 40 0, clip]{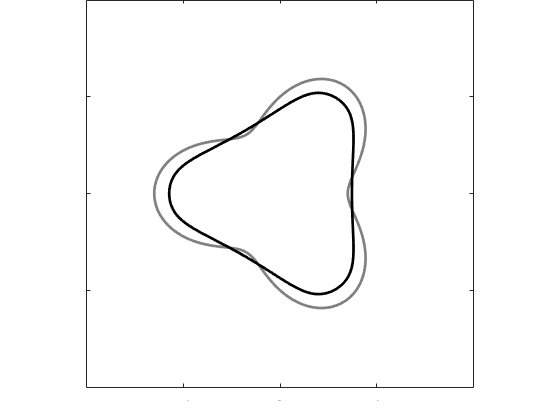}
\end{subfigure}
\begin{subfigure}[b]{0.3\textwidth}
\includegraphics[width=\textwidth,trim=40 20 40 0, clip]{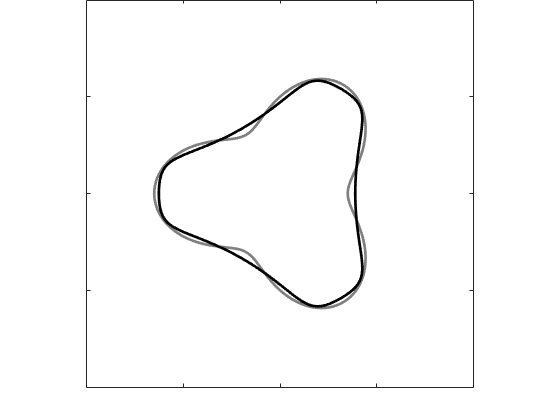}
\end{subfigure}
\caption{The initial  guess for sinusoidal perturbation of a disk. The black curve is the reconstructed curve and the gray curve is the actual shape. Starting from the left side to the right side, the conductivity is 3, 10, 50, and the first row shows the equivalent ellipse and the second row shows the reference shape.}
\label{fig:sinu_initial}
\end{figure}

\begin{figure}[H]	
\begin{subfigure}[b]{0.3\textwidth}
\includegraphics[width=\textwidth,trim=40 20 40 0, clip]{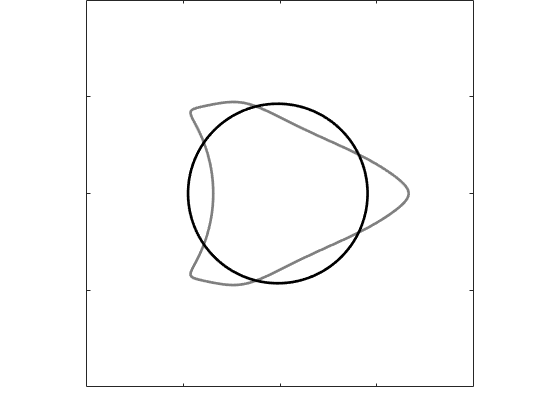}
\end{subfigure}
\begin{subfigure}[b]{0.3\textwidth}
\includegraphics[width=\textwidth,trim=40 20 40 0, clip]{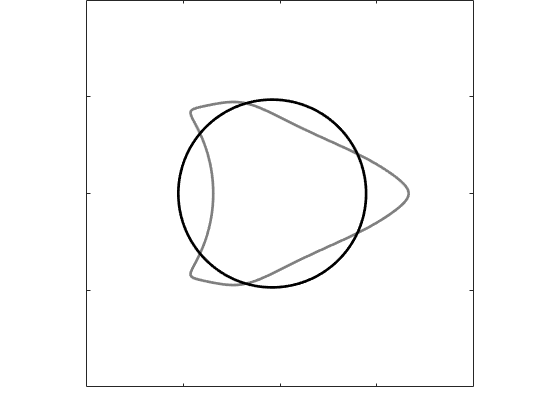}
\end{subfigure}
\begin{subfigure}[b]{0.3\textwidth}
\includegraphics[width=\textwidth,trim=40 20 40 0, clip]{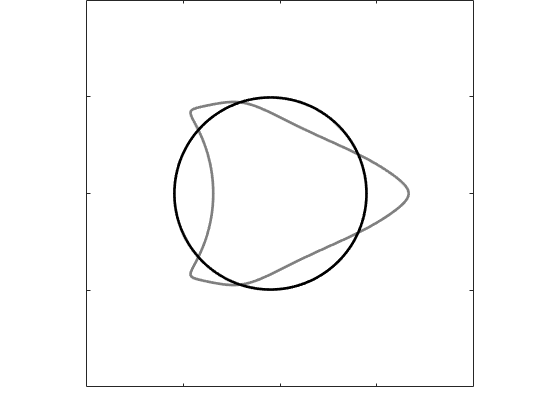}
\end{subfigure}
\vskip 0.2cm
\begin{subfigure}[b]{0.3\textwidth}
\includegraphics[width=\textwidth,trim=40 20 40 0, clip]{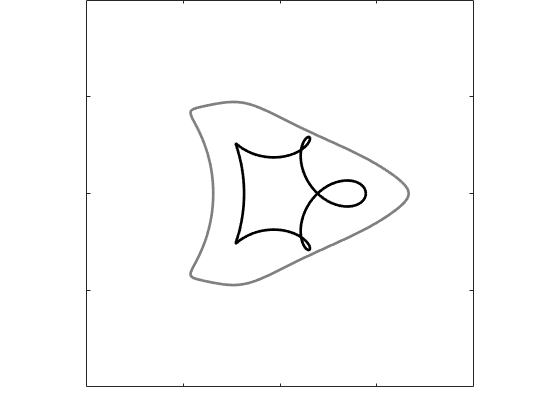}
\end{subfigure}
\begin{subfigure}[b]{0.3\textwidth}
\includegraphics[width=\textwidth,trim=40 20 40 0, clip]{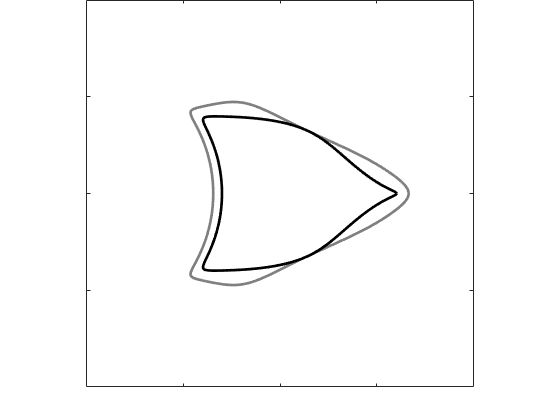}
\end{subfigure}
\begin{subfigure}[b]{0.3\textwidth}
\includegraphics[width=\textwidth,trim=40 20 40 0, clip]{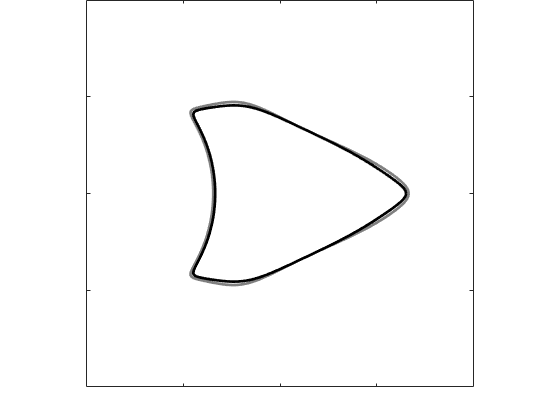}
\end{subfigure}
\caption{The initial  guess for the kite-shaped domain. The black curve is the reconstructed curve and the gray curve is the actual shape. Starting from the left side to the right side, the conductivity is 0.5, 0.1, 0.02, and the first row shows the equivalent ellipse and the second row shows the reference shape.}
\label{fig:kite_initial}
\end{figure}

Let's compare two methods to make the initial guess. Both methods have advantages and disadvantages. Using the equivalent ellipse is always stable, but it does not contain information of high order GPTs. Using approximation from FPTs is almost exact for extreme conductivities, but the error increases as conductivity closes to 1. 
Figure \ref{fig:sinu_initial} shows the initial guess for the perturbed circle with perturbation $0.3\cos 3\theta$ on the radius $r$ and GPTs of order 6 are given. When conductivity is 3, the equivalent ellipse is more accurate. But as conductivity goes to infinite, the initial guess approximating from FPTs becomes similar to the target domain. The target domain of Figure \ref{fig:kite_initial} is the kite-shaped domain and GPTs of order 6 are given. When the conductivity is 0.5, the initial guess approximating from FPTs is not even simply connected. But as conductivity goes to zero, the initial guess approximating from FPTs becomes almost same as the target domain.

\subsection{Recursive scheme}\label{sec:recursion}

We now provide the optimization method which uses the cost function defined by the GPTs. 
We let $B$ denote the target domain.

For the optimization process, we use the cost function $J_K[D]$ given by
$$J_K[D] := \frac{1}{2} \sum_{|\alpha|, |\beta| \le K} \Bigg|\sum_{\alpha, \beta} a_{\alpha} b_{\beta} M_{\alpha \beta}  (D,\lambda) - \sum_{\alpha, \beta} a_{\alpha} b_{\beta} M_{\alpha \beta} (B,\lambda) \Bigg|^2, $$
where $a_\alpha$ and $b_\beta$ are coefficients of the harmonic polynomials $H = \sum_\alpha a_\alpha x^\alpha$ and $F = \sum_\beta b_\beta x^\beta$. Note that the cost function vanishes only when each $M_{\alpha \beta}$ is same. To minimize the cost function, we use following lemma.

\begin{lemma}[\cite{ammari:2012:GPT}]\label{lemma:shapederivative}
Let $\p D_\ep(x) := \{ x + \ep h(x) \nu(x) \, | \,  x \in \p D\}$.  Suppose that $a_\alpha$ and $b_\beta$ are constants such that $H=\sum_\alpha a_\alpha x^\alpha$ and $F=\sum_\beta b_\beta x^\beta$ are harmonic polynomials. Then
	\begin{align*}
	&\sum_{\alpha, \beta} a_\alpha b_\beta M_{\alpha \beta} (D_\ep,\sigma_0) - \sum_{\alpha, \beta} a_\alpha b_\beta M_{\alpha \beta} (D,\sigma_0) 	\\
=& \ep (\sigma_0-1) \int_{\p D} h(x) \left[ \frac{\p v}{\p \nu} \bigg|^- \frac{\p u}{\p \nu} \bigg|^- + \frac{1}{\sigma_0} \frac{\p v}{\p T} \bigg|^- \frac{\p u}{\p T} \bigg|^- \right] (x) d \sigma(x) + O(\ep^2),
	\end{align*}
	where $u$ and $v$ are solutions to the following transmission problems:
	\beq\label{cond_eqn1}
	\begin{cases}
		\ds\Delta u=0\quad&\mbox{in } D\cup (\RR^2\setminus\overline{D}), \\
		\ds u\big|^+=u\big|^-\quad&\mbox{on }\p D,\\
		\ds \pd{u}{\nu}\bigg|^+=\sigma_0\pd{u}{\nu}\bigg|^-,\quad&\mbox{on }\p D,\\
		\ds u(x) - H(x)  =O({|x|^{-1}})\quad&\mbox{as } |x| \to \infty
	\end{cases}
\end{equation}
and
\beq\label{cond_eqn2}
\begin{cases}
	\ds\Delta v=0\quad&\mbox{in } D\cup (\RR^2\setminus\overline{D}), \\
	\ds \sigma_0v\big|^+=v\big|^-\quad&\mbox{on }\p D,\\
	\ds \pd{v}{\nu}\bigg|^+=\pd{v}{\nu}\bigg|^-,\quad&\mbox{on }\p D,\\
	\ds v(x) - F(x)  =O({|x|^{-1}})\quad&\mbox{as } |x| \to \infty,
\end{cases}
\end{equation}
where $H(x)=\sum a_\alpha x^\alpha$ and $F(x)=\sum b_\beta x^\beta$ are harmonic polynomials.
\end{lemma} 

From Lemma \ref{lemma:shapederivative}, the shape derivative of $J_K [D]$ becomes
$$
\la d_S J_K [D],h\ra = \sum_{|\alpha|, |\beta| \le K} \delta_{HF}\la \phi_{HF} ,h\ra_{L^2(\p D)},
$$
with
\begin{align*}
\delta_{HF} &= \sum_{\alpha, \beta} a_{\alpha} b_{\beta} M_{\alpha \beta}  (D,\lambda) - \sum_{\alpha, \beta} a_{\alpha} b_{\beta} M_{\alpha \beta} (B,\lambda), \\
\phi_{HF} &= (\sigma_0-1)\left[\pd{u}{\nu}\bigg|^-\pd{v}{\nu}\bigg|^-+\frac{1}{\sigma_0}\pd{u}{T}\bigg|^-\pd{v}{T}\bigg|^-\right](x),\quad x\in\p D.
\end{align*}
Here, $u$ and $v$ satisfy \eqnref{cond_eqn1} and \eqnref{cond_eqn2} for $H = \sum_\alpha a_\alpha x^\alpha$ and $F = \sum_\beta b_\beta x^\beta$.

For each step, we modify the shape using the gradient descent method with the cost function when $\widetilde{D} = D_{init}$. The formula becomes
$$\p D_{mod} = \p D_{init} - \Bigg( \frac{J_K [D_{init}]  }{\sum_j \la d_s J_K [D_{init}], \varphi_j \ra^2} \sum_j \la d_s J_K [D_{init}], \varphi_j \ra \varphi_j \Bigg)\nu, $$
where $\{\varphi_j\}$ is the basis on $L^2(\p D_{init})$. From Lemma \ref{lemma:shapederivative}, we only know the information about $h$ is the inner product with $\phi_{HF}$. Hence, we take the basis set $\{\varphi_i\}$ as $\phi_{HF}$'s  for harmonic polynomials $H = \sum_\alpha a_\alpha x^\alpha$ and $F = \sum_\beta b_\beta x^\beta$.

\subsection{Numerical results}

\begin{figure}[H]
\begin{subfigure}[b]{0.3\textwidth}
\includegraphics[width=\textwidth,trim=40 20 40 0, clip]{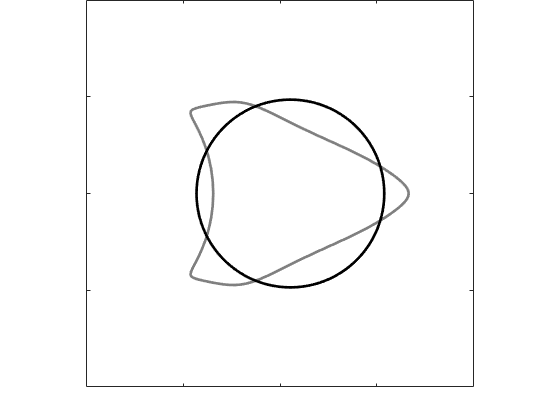}
\end{subfigure}
\hskip 1cm
\begin{subfigure}[b]{0.3\textwidth}
\includegraphics[width=\textwidth,trim=40 20 40 0, clip]{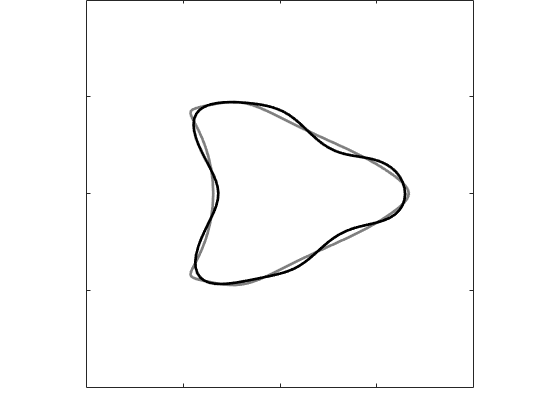}
\end{subfigure}
\vskip 0.2cm
\begin{subfigure}[b]{0.3\textwidth}
\includegraphics[width=\textwidth,trim=40 20 40 0, clip]{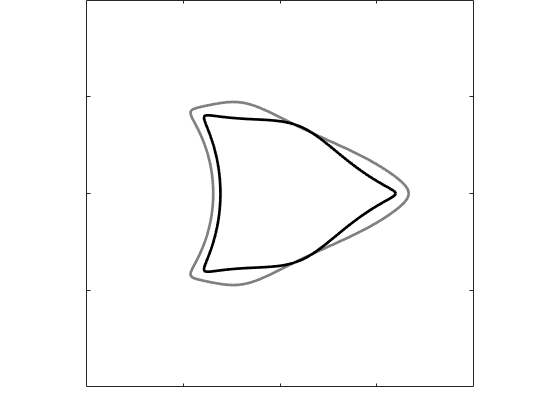}
\end{subfigure}
\hskip 1cm
\begin{subfigure}[b]{0.3\textwidth}
\includegraphics[width=\textwidth,trim=40 20 40 0, clip]{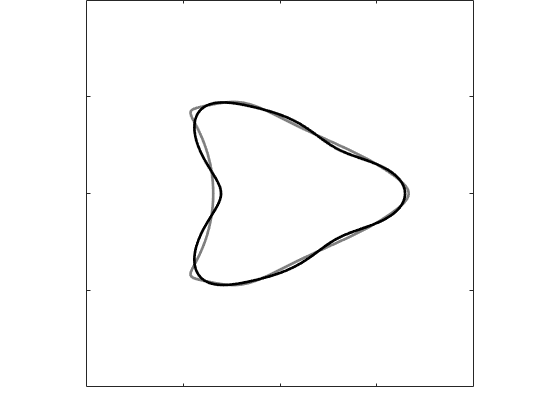}
\end{subfigure}
\caption{The reconstructed inclusion for the kite-shaped target with conductivity 10. The left figures are the initial guesses and the right figures are the reconstructed results. The first row is shape reconstruction starting from the equivalent ellipse and the second row is reconstruction approximating from FPTs.}
\label{recon:kite_cond10}
\end{figure}

The Figure \ref{recon:kite_cond10} shows the reconstructed results for kite-shaped domain with the initial guess given by the reference shape. 

\begin{figure}[h]
	\begin{subfigure}[b]{0.3\textwidth}
		\includegraphics[width=\textwidth,trim=40 20 40 0, clip]{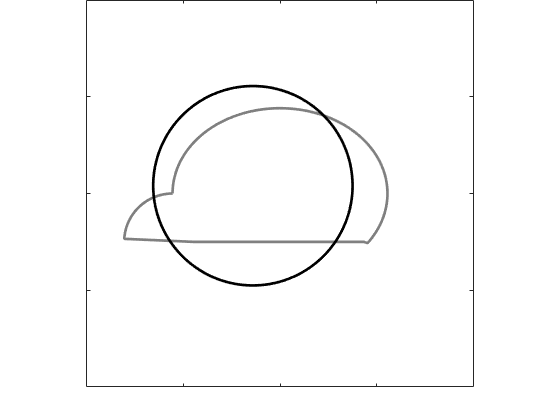}
	\end{subfigure}
	\hskip 1cm
	\begin{subfigure}[b]{0.3\textwidth}
		\includegraphics[width=\textwidth,trim=40 20 40 0, clip]{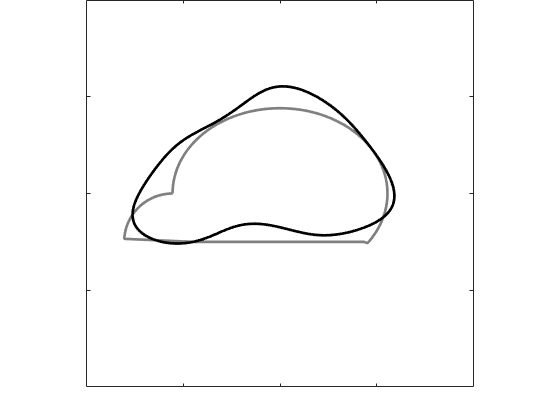}
	\end{subfigure}
	\vskip 0.2cm
	\begin{subfigure}[b]{0.3\textwidth}
		\includegraphics[width=\textwidth,trim=40 20 40 0, clip]{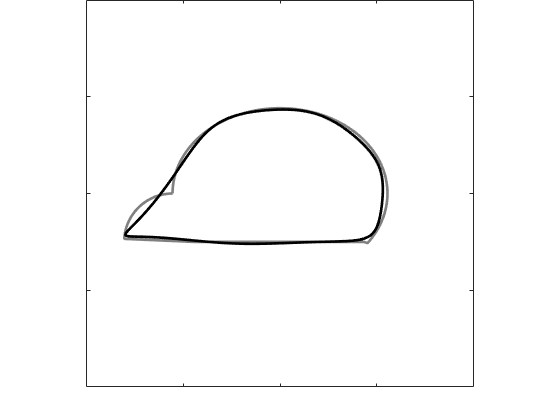}
	\end{subfigure}
	\hskip 1cm
	\begin{subfigure}[b]{0.3\textwidth}
		\includegraphics[width=\textwidth,trim=40 20 40 0, clip]{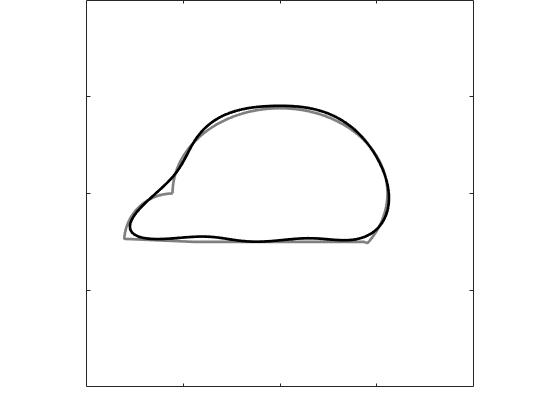}
	\end{subfigure}
	\caption{The reconstructed inclusion for the cap-shaped target with conductivity 50. The left figures are the initial guesses and the right figures are the reconstructed results. The first row is shape reconstruction starting from equivalent ellipse and the second row is reconstruction starting from approximation from FPTs.}
	\label{recon:cap_cond50}
\end{figure}

The Figure \ref{recon:cap_cond50} reveals the reconstructed results for the cap-shaped domain. Figure \ref{recon:cap_cond50} shows the importance of the initial guess.


\begin{thebibliography}{10}

\bibitem{ammari:2013:MSM}
Habib Ammari, Josselin Garnier, Wenjia Jing, Hyeonbae Kang, Mikyoung Lim, Knut
  S{\o}lna, and Han Wang.
\newblock {\em Mathematical and Statistical Methods for Multistatic Imaging},
  volume 2098 of {\em Lecture Notes in Mathematics}.
\newblock Springer, Cham, 2013.

\bibitem{ammari:2014:GPT}
Habib Ammari, Josselin Garnier, Hyeonbae Kang, Mikyoung Lim, and Sanghyeon Yu.
\newblock Generalized polarization tensors for shape description.
\newblock {\em Numerische Mathematik}, 126(2):199--224, Feb 2014.

\bibitem{ammari:2004:RSI}
Habib Ammari and Hyeonbae Kang.
\newblock {\em Reconstruction of small inhomogeneities from boundary
  measurements}, volume 1846 of {\em Lecture Notes in Mathematics}.
\newblock Springer, Berlin, Heidelberg, 2004.

\bibitem{ammari:2007:PMT}
Habib Ammari and Hyeonbae Kang.
\newblock {\em Polarization and Moment Tensors: With Applications to Inverse
  Problems and Effective Medium Theory}, volume 162 of {\em Applied
  Mathematical Sciences}.
\newblock Springer Science+Business Media, 2007.

\bibitem{ammari:2012:GPT}
Habib Ammari, Hyeonbae Kang, Mikyoung Lim, and Habib Zribi.
\newblock The generalized polarization tensors for resolved imaging. part i:
  Shape reconstruction of a conductivity inclusion.
\newblock {\em Mathematics of Computation}, 81(277):367--386, 2012.

\bibitem{choi:2018:CEP}
D.~Choi, J.~Helsing, and M.~Lim.
\newblock Corner effects on the perturbation of an electric potential.
\newblock {\em SIAM Journal on Applied Mathematics}, 78(3):1577--1601, 2018.

\bibitem{choi:2018:GME}
Doosung Choi, Junbeom Kim, and Mikyoung Lim.
\newblock Geometric multipole expansion and its application to neutral
  inclusions of general shape.
\newblock {\em arXiv preprint arXiv:1808.02446}, 2018.

\bibitem{duren:1983:UF}
Peter~L. Duren.
\newblock {\em Univalent Functions}, volume 259 of {\em Grundlehren der
  mathematischen Wissenschaften}.
\newblock Springer-Verlag New York, 1983.

\bibitem{escauriaza:1992:RTW}
Luis Escauriaza, Eugene~B. Fabes, and Gregory Verchota.
\newblock On a regularity theorem for weak solutions to transmission problems
  with internal lipschitz boundaries.
\newblock {\em Proceedings of the American Mathematical Society},
  115(4):1069--1076, 1992.

\bibitem{escauriaza:1993:RPS}
Luis Escauriaza and Jin~Keun Seo.
\newblock Regularity properties of solutions to transmission problems.
\newblock {\em Transactions of the American Mathematical Society},
  338(1):405--430, 1993.

\bibitem{faber:1903:UPE}
Georg Faber.
\newblock Über polynomische entwickelungen.
\newblock {\em Mathematische Annalen}, 57:389--408, 1903.

\bibitem{helsing:2013:SIE}
Johan Helsing.
\newblock Solving integral equations on piecewise smooth boundaries using the
  rcip method: A tutorial.
\newblock {\em Abstract and Applied Analysis}, vol. 2013, Article ID 938167:20
  pages, 2013.

\bibitem{helsing:2017:CSN}
Johan Helsing, Hyeonbae Kang, and Mikyoung Lim.
\newblock Classification of spectra of the neumann–poincaré operator on
  planar domains with corners by resonance.
\newblock {\em Annales de l'Institut Henri Poincare (C) Non Linear Analysis},
  34(4):991 -- 1011, 2017.

\bibitem{jung:2018:NSS}
YoungHoon Jung and Mikyoung Lim.
\newblock A new series solution method for the transmission problem.
\newblock {\em arXiv preprint arXiv:1803.09458}, 2018.

\bibitem{kellogg:2012:FPT}
Oliver~Dimon Kellogg.
\newblock {\em Foundations of Potential Theory}, volume~31 of {\em Die
  Grundlehren der Mathematischen Wissenschaften}.
\newblock Springer, Berlin, Heidelberg, 1967.

\bibitem{verchota:1984:LPR}
Gregory Verchota.
\newblock Layer potentials and regularity for the dirichlet problem for
  laplace's equation in lipschitz domains.
\newblock {\em Journal of Functional Analysis}, 59(3):572 -- 611, 1984.

\end{thebibliography}
\end{document}